\documentclass[12 pt, draftcls, onecolumn]{IEEEtran}
\usepackage{amssymb,epsfig,color,cite,amsmath,amsfonts,mathrsfs,algorithm,algorithmic}
\usepackage{epstopdf}

\newtheorem{lemma}{Lemma}[section]
\newtheorem{theorem}{Theorem}[section]

\newtheorem{remark}{Remark}[section]

\newtheorem{assumption}{Assumption}[section]
\allowdisplaybreaks

\begin{document}
%
\title{\huge Distributed Continuous-Time Algorithm for Constrained Convex Optimizations via Nonsmooth Analysis Approach}
\author{Xianlin~Zeng, Peng~Yi, and Yiguang~Hong
\thanks{This work was supported by NSFC (61603378, 61573344), Beijing Natural Science Foundation (4152057), Project 2014CB845301/2, and China Postdoctoral Science Foundation (2015M581190).}
\thanks{X. Zeng and Y. Hong are with Key Laboratory of Systems and Control, Institute of Systems Science, Chinese Academy of Sciences, Beijing, 100190, China (e-mail: xianlin.zeng@amss.ac.cn; yghong@iss.ac.cn).}
\thanks{P.~Yi  is with the Department of Electrical and Computer Engineering, University of Toronto,  Toronto,   Canada  (e-mail: yi.peng@scg.utoronto.ca).}
}

\markboth{Submitted to IEEE Transactions on Automatic Control}{}

\maketitle

\begin{abstract}
This technical note studies the distributed optimization problem of a sum of nonsmooth convex cost functions with local constraints. At first, we propose a novel distributed continuous-time projected algorithm, in which each agent knows its local cost function and local constraint set, for the constrained optimization problem.
Then we prove that all the agents of the algorithm can find the same optimal solution,
and meanwhile, keep the states bounded while seeking the optimal solutions.
We conduct a complete convergence analysis by employing nonsmooth Lyapunov functions for the stability analysis of differential inclusions.
Finally, we provide a numerical example for illustration.
\end{abstract}


%

\section{Introduction}
%
%
%
%
The distributed optimization of a sum of convex functions is an important class of decision and data processing problems over network systems,
and has been intensively studied in recent years (see \cite{NOP:2010,LO:2011,WE:2011,GC:2014,SJH:TAC:2013,YHL:SCL:2015} and references therein).
In addition to the discrete-time distributed optimization algorithms (e.g., \cite{NOP:2010,LO:2011}), 
continuous-time multi-agent solvers have recently been applied to distributed optimization problems as a promising and useful technique \cite{WE:2011,SJH:TAC:2013,GC:2014,KCM:A:2015,YHL:SCL:2015,QLX:A:2016}, thanks to the well-developed continuous-time stability theory.

Constrained distributed optimization, in which the feasible solutions are limited to a
certain region or range, is significant in a number of network decision applications, including multi-robot motion planning, resource allocation in communication networks, and economic dispatch in power grids.
{In practice, local constraints in the distributed optimization design are often necessary due to the performance limitations of the agents in computation and communication capacities as well as task requirements of privacy and security. For example, in large-scale optimization problems, the computation/communication capacity of a single agent may not be enough to handle all the constraints of the agents;
in alignment or resource allocation problems, each agent's feasible choice is limited to a certain range, while the agents may not want to share their private information with others; and
in strategic social networks, the agents keep their own limit constraints or budget constraints confidential for security concerns.}
However, due to the consideration of local constraints, the design of such algorithms, to minimize the global cost functions within the feasible set while allowing the agents operate with only local cost functions and local constraints, 
is a very difficult task. Conventionally, the projection method has been widely adopted in the algorithm design for constrained optimization \cite{AHU:1972,Ruszczynski:2006} and related problems \cite{NZ:1996}.
\cite{YHL:SCL:2015} constructed a primal-dual type continuous-time projected algorithm to solve a distributed optimization problem, where each agent has its own private constraint function, while \cite{QLX:A:2016} proposed a continuous-time distributed projected dynamics for constrained optimization, where the agents share the same constraint set.  Moreover,
\cite{LW:TAC:2015} presented a primal-dual continuous-time projected algorithm for distributed nonsmooth optimization, where each agent has its own local bounded constraint set,
though its auxiliary variables may be asymptotically unbounded.

The purpose of this technical note is to propose a novel continuous-time projected algorithm for distributed nonsmooth convex optimization problems where each agent has its own general local constraint set. The main contributions of the note are four folds.
Firstly, a distributed continuous-time algorithm is proposed for the agents to find the same optimal solution based only on local cost functions and local constraint sets, by combining primal-dual methods for saddle point seeking and projection methods for set constraints.  The proposed algorithm is consistent with those in \cite{GC:2014,WE:2011,SJH:TAC:2013} when there were no constraints in the optimization problem.
Secondly, nonsmooth cost functions are considered here, while smooth cost functions were discussed in most continuous-time distributed optimization designs \cite{YHL:SCL:2015,KCM:A:2015}.  To solve the complicated problem, nonsmooth Lyapunov functions are employed along with the stability theory of differential inclusions (resulting from the nonsmooth cost functions) to conduct a complete and original convergence analysis.
Thirdly, our proposed algorithm is proved to solve the optimization problem and have bounded states while seeking the optimal solutions, and therefore, further improves the recent interesting result in \cite{LW:TAC:2015}, whose algorithm may have asymptotically unbounded states.
Finally, different from the strict/strong convexity in existing results \cite{YHL:SCL:2015,KCM:A:2015}, general convexity is investigated.  In fact, our nonsmooth analysis techniques also guarantee the convergence of the algorithm  even when the problem has a continuum of optimal solutions due to the convexity. Therefore, the convergence analysis provides additional insights and understandings for continuous-time distributed optimization algorithms compared with \cite{WE:2011,SJH:TAC:2013,YHL:SCL:2015,KCM:A:2015}.

The remainder of this note is organized as follows. In Section \ref{sec:def}, notations and definitions are presented and reviewed.
In Section \ref{Distributed_Optimization}, a constrained convex (nonsmooth) optimization problem is formulated and a distributed continuous-time projected algorithm is proposed. In Section \ref{convergence_simulation}, a complete proof is presented to show that the algorithm state is bounded and the agents' estimates are convergent to the same optimal solution, and simulation studies are carried out for illustration. Finally, in Section \ref{conclusion}, concluding remarks are given.

\section{Mathematical Preliminaries}
\label{sec:def}
In this section, we introduce necessary notations, definitions and preliminaries about graph theory
and projection operators.

\subsection{Notations}
Let $\mathbb R$ denote the set of real numbers; let $\mathbb R^n$ and $\mathbb R^{n\times m}$
denote the set of $n$-dimensional real column vectors and the set of $n$-by-$m$ real matrices, respectively; $\mathfrak B(\mathbb R^q)$ denotes the collection of all subsets of $\mathbb R^q$; $I_n$ denotes the $n\times n$ identity matrix and $(\cdot)^\mathrm{T}$ denotes the transpose.  Furthermore, $\|\cdot\|$ denotes the Euclidean norm.
Write $\mathrm {rank}(A)$ for the rank of a matrix $A$, $\mathrm{range} (A)$ for the range of  $A$, $\ker (A)$ for the kernel of  $A$, $\lambda _{\mathrm {max}} (A)$ for the largest eigenvalue of  $A$, $\mathbf 1_n$ for the $n\times 1$ ones vector, $\mathbf 0_n$ for the $n\times 1$ zeros vector, and $A \otimes B$ for the Kronecker product of matrices $A$ and $B$.
Denote
$A> 0$ (or $A\geq 0$) when matrix $A\in\mathbb R^{n\times n}$ is positive definite (or positive semi-definite).  Also,
denote $\overline{\mathcal{S}}$ as the closure of a subset $\mathcal{S}\subset\mathbb{R}^{n}$, $\mathrm{int}(\mathcal{S})$ as the interior of  $\mathcal{S}$,
$\mathcal N_{\mathcal S}({x})$ as the normal cone of $\mathcal S$ at an element ${x}\in\mathcal S$, $\mathcal T_{\mathcal S}({x})$ as the tangent cone of $\mathcal S$ at an element ${x}\in\mathcal S$, and
$\mathcal B_{\epsilon}(p),p\in\mathbb R^n$ as the open ball {\em centered} at $p$ with {\em radius} $\epsilon>0$. Denote $\mathrm{dist}(x,\mathcal M)$ as the distance from a point $x$ to a set $\mathcal M$ (that is, $\mathrm{dist}(x,\mathcal M) \triangleq \inf_{p\in\mathcal M}\|p-x\|$), and $x(t)$ approaches $\mathcal M$ if $x(t)\rightarrow \mathcal M$ as $t\rightarrow \infty$
(that is, for each $\epsilon>0$, there is $T>0$ such that $\mathrm{dist}(x(t),\mathcal M)<\epsilon$ for all $t>T$).

\subsection{Graph Theory}
A weighted undirected graph $\mathcal G$  is denoted by $\mathcal G(\mathcal V, \mathcal E, A)$, where $\mathcal V=\{1,\ldots, n\}$ is a set of nodes, $\mathcal E\subset\mathcal V \times \mathcal V$ is a set of edges, $ A=[a_{i,j}]\in\mathbb R^{n\times n}$ is a {\em {weighted} adjacency matrix} such that {$a_{i,j},\,a_{j,i}>0$} if ${\{i,j\}}\in\mathcal E,j\not=i$, and $a_{i,j}=0$ otherwise. 
The {\em {weighted} Laplacian matrix} is $L_n=D- A$, where $D\in\mathbb R^{n\times n}$ is diagonal with $D_{i,i}=\sum^n_{j=1,j\not=i} a_{i,j}$, $i\in\{1,\ldots,n\}$. {In this note, we  call $L_n$ the Laplacian matrix and $A$ the adjacency matrix of $\mathcal G$ for convenience when there is no confusion.}
Specifically, if the weighted undirected graph $\mathcal G$ is connected,
then $L_n\geq 0$, ${\mathrm{rank} (L_n)}=n-1$, and $\ker (L_n)=\{k\mathbf{1}_n:k\in\mathbb R\}$.

\subsection{Projection Operator}

Define $P_{K}(\cdot)$ as a projection operator given by $P_{K}(u)=\mathrm{arg}\,\min_{v\in K}\|u-v\|$, where $K\subset \mathbb R^n$.
\begin{lemma}\cite{KS:1982}
If $K\subset \mathbb R^n$ is a closed convex set, then
\begin{eqnarray}\label{projection_inequality}
(u-P_{K}(u))^{\rm T}(v-P_{K}(u))\leq 0,\quad \forall u\in\mathbb R^n,\quad \forall v\in K.
\end{eqnarray}
\end{lemma}

\section{Problem Description and Optimization Algorithm}
\label{Distributed_Optimization}
\subsection{Problem Description}
Consider a network of $n$ agents interacting over a graph $\mathcal G$.
There is a local cost function $f^{i}:\mathbb R^q\rightarrow\mathbb R$ and a local feasible constraint set $\Omega_i\subset\mathbb R^q$ for all $i\in\{1,\ldots,n\}$.
The global cost function of the network is $f(x) = \sum_{i=1}^n f^i(x)$, and the feasible set is the intersection of local constraint sets, that is, $x\in \Omega_0\triangleq\bigcap_{i=1}^n\Omega_i\subset\mathbb R^q$.
Then we will provide a distributed algorithm to solve
\begin{eqnarray}\label{optimization_problem}
\min_{x\in\Omega_0}f(x), \quad f(x)= \sum_{i=1}^n f^i(x),\quad x\in\Omega_0\subset\mathbb R^q,
\end{eqnarray}
where each agent only uses  its own local cost function, its local constraint, and the shared information of its neighbors through {constant local communications}.

To ensure the wellposedness of problem (\ref{optimization_problem}), the following assumption is needed.
\begin{assumption}\label{Assumption1}~

\begin{enumerate}
  \item The weighted graph $\mathcal G$ is connected and undirected with symmetric weighted Laplacian matrix $L_n$.
  \item For all $i\in\{1,\dots,n\}$, $f^i$ is continuous and convex {on an open set containing $\Omega_i$}, and $\Omega_i\subset\mathbb R^q$ is {closed} and {convex} with $\bigcap_{i=1}^n \mathrm{int}(\Omega_i)\not=\emptyset$.
  \item There exists at least one  optimal solution to problem (\ref{optimization_problem}).
\end{enumerate}
\end{assumption}

\begin{remark}
Problem (\ref{optimization_problem}) covers many problems in  recent distributed optimization studies. For example, it introduces the constraints compared with the unconstrained optimization model in \cite{GC:2014}.
Moreover, it generalizes the model in \cite{QLX:A:2016} by allowing heterogeneous constraints, and extends the models in \cite{YHL:SCL:2015} and \cite{LW:TAC:2015}, which considered function constraints and hyper box (sphere) constraints, respectively.
\hfill$\blacklozenge$
\end{remark}

Let $x_i(t)\in\Omega_i\subset\mathbb R^q$ be the estimate of agent $i$ at time instant $t\geq 0$ for the optimal solution. Let $\mathbf{L}\triangleq L_n\otimes I_q\in\mathbb R^{nq\times nq}$, where $L_n\in\mathbb R^{n\times n}$ is the
Laplacian matrix of $\mathcal G$. Denote $\mathbf{x} \triangleq [x_1^{\rm T},\ldots,x_n^{\rm T}]^{\rm T}\in\Omega\subset\mathbb R^{nq}$ and $\mathbf{f}(\mathbf{x})\triangleq\sum_{i=1}^n f^i(x_i)$ with $\mathbf{x}\in\Omega$, where $\Omega\triangleq\prod_{i=1}^n \Omega_i$ is the Cartesian product of $\Omega_i,\,i\in\{1,\ldots,n\}$. Then we arrive at the following lemma by directly analyzing the optimality condition.

\begin{lemma}\label{optimal_eqv}
Suppose Assumption \ref{Assumption1} holds and $\alpha>0$. ${x}^*\in\Omega_0\subset\mathbb R^q$ is an optimal solution to problem (\ref{optimization_problem})
if and only if there exist $\mathbf{x}^*=\mathbf{1}_n \otimes x^*\in\Omega\subset\mathbb{R}^{nq}$
and $\lambda^*\in\mathbb{R}^{nq}$ such that
\begin{subequations}\label{optimal_condition}
\begin{align}
&\mathbf{0}_{nq}\in  \big{\{}P_{\mathcal T_{\Omega} (\mathbf{x}^*)}(-g(\mathbf{x}^*)-\alpha\mathbf{L\lambda}^*):
                    \,g(\mathbf{x}^*)\in\partial\mathbf{f}(\mathbf{x}^*)\big{\}},\label{optimal_condition_1}\\
&\mathbf{Lx}^* = \mathbf{0}_{nq},\label{optimal_condition_2}
\end{align}
\end{subequations}
where $\mathcal T_{\Omega}(\mathbf{x}^*)$ is the tangent cone of $\Omega$ at an element $\mathbf{x}^*\in\Omega$ and $P_{\mathcal T_{\Omega}(\mathbf{x}^*)}(\cdot)$ is the projection operator to $\mathcal T_{\Omega}(\mathbf{x}^*)$.
\end{lemma}

\begin{IEEEproof}
According to Theorem 3.33 in \cite{Ruszczynski:2006}, $x^*$ is an optimal solution to problem \eqref{optimization_problem} if and only if
\begin{equation}
\mathbf{0}_q \in  \partial f(x^*) + \mathcal N_{\Omega_0}(x^*),\label{optimal_orginal}
\end{equation}
where $\mathcal N_{\Omega_0}({x}^*)$ is the normal cone of $\Omega_0$ at ${x}^*\in\Omega_0=\bigcap_{i=1}^n\Omega_i$.
Note that $f^i(\cdot),\,i=1,\ldots,n,$ is convex and $\bigcap_{i=1}^n \mathrm{int}(\Omega_i)\not=\emptyset$ by Assumption \ref{Assumption1}. It follows from Theorem 2.85 and Lemma 2.40 in \cite{Ruszczynski:2006} that $\partial f(x^*)= \sum_{i=1}^n \partial f^i(x^*)$ and $\mathcal N_{\Omega_0}(x^*)= \sum_{i=1}^n \mathcal N_{\Omega_i}(x^*)$. To prove this lemma, one only needs to show \eqref{optimal_orginal} holds if and only if (\ref{optimal_condition}) is satisfied.

Suppose (\ref{optimal_condition}) holds. Since graph $\mathcal G$ is connected,  there exists $x^*\in\mathbb R^q$ such that $\mathbf{x}^*=\mathbf{1}_n \otimes x^*\in\mathbb R^{nq}$ because of (\ref{optimal_condition_2}). Note that $\mathbf{0}_{nq}=P_{\mathcal T_{\Omega} (\mathbf{x}^*)}(-g(\mathbf{x}^*)-\alpha\mathbf{L\lambda}^*)$ if and only if $-g(\mathbf{x}^*)-\alpha\mathbf{L\lambda}^*\in\mathcal N_{\Omega}(\mathbf{x}^*)$.
Let $a_{i,j}$ be the $(i,j)$th entry of the
adjacency matrix of $\mathcal G$ and $ \lambda^* = [(\lambda_1^*)^{\rm T},\ldots,(\lambda_n^*)^{\rm T}]^{\rm T}\in\mathbb R^{nq}$ with $\lambda^*_i\in\mathbb R^q,\,i\in\{1,\ldots,n\}$.
Then $\eqref{optimal_condition_1}$ holds if and only if there exists $g_i(x^*)\in \partial f^i(x^*)$ such that
$-g_i(x^*)-\alpha \sum_{j=1i}^n a_{i,j}(\lambda_i^*-\lambda_j^*) \in\mathcal N_{\Omega_i}(x^*), i=1,...,n$.
Because $L_n=L_n^{\rm T}$ {by Assumption \ref{Assumption1}}, $\sum_{i=1}^n\sum_{j=1}^n a_{i,j}(\lambda_i^*-\lambda_j^*)=1/2\sum_{i=1}^n\sum_{j=1}^n (a_{i,j}-a_{j,i})(\lambda_i^*-\lambda_j^*)=\mathbf{0}_q$ 
and  $-\sum_{i=1}^n g_i(x^*) \in \sum_{i=1}^n\mathcal N_{\Omega_i}(x^*) =\mathcal N_{\Omega_0}(x^*) $. Since $\sum_{i=1}^n g_i(x^*)\in\sum_{i=1}^n \partial f^i(x^*)=\partial f(x^*)$, \eqref{optimal_orginal} is thus proved.

Conversely, suppose (\ref{optimal_orginal}) holds. Let $\mathbf{x}^*=\mathbf{1}_n \otimes x^*$. (\ref{optimal_condition_2}) is clearly true. It follows from (\ref{optimal_orginal}) that there exists $g_i(x^*)\in \partial f^i(x^*)$ such that $  -\sum_{i=1}^n g_i(x^*) \in \sum_{i=1}^n \mathcal{N}_{\Omega_i}(x^*)$. Choose $z_i(x^*) \in \mathcal{N}_{\Omega_i}(x^*) ,\,i=1,\ldots,n$, such that $-\sum_{i=1}^n  g_i(x^*)  =  \sum_{i=1}^n z_i(x^*)  $.
Next, define vectors  $l_i(x^*)\triangleq z_i(x^*)+g_i(x^*),\,i=1,\ldots,n$.
It is clear that $\sum_{i=1}^n l_i(x^*) =\mathbf{0}_q $.
{Note that $\mathbf{L}$ is symmetric by Assumption \ref{Assumption1}.}
By the fundamental theorem of linear algebra, the sets $\mathrm{ker}(\mathbf{L})$ and $\mathrm{range}(\mathbf{L})$ form an orthogonal decomposition of $\mathbb R^{nq}$. Define $l(x^*)\triangleq [l_1(x^*)^{\rm T},...,l_n(x^*)^{\rm T}]^{\rm T}\in\mathbb R^{nq}$. For all $\mathbf{x}=\mathbf{1}_n \otimes x\in\ker(\mathbf{L})$, $l(x^*)^{\rm T}\mathbf{x}=\sum_{i=1}^n l_i(x^*)^{\rm T}x=0$ and hence, $l(x^*)\in \mathrm{range}(\mathbf{L})$ and there exists {$\lambda^*\in\mathbb R^{nq}$} such that $l(x^*)=-\alpha\mathbf{L}\lambda^{*}$.
Thus, there exists $ \lambda^* = [(\lambda_1^*)^{\rm T},\ldots,(\lambda_n^*)^{\rm T}]^{\rm T}\in\mathbb R^{nq}$ with $\lambda^*_i\in\mathbb R^q$ such that $-g_i(x^*)-\alpha \sum_{j=1}^n a_{i,j}(\lambda_i^*-\lambda_j^*)=-g_i(x^*)+l_i(x^*)=z_i(x^*) \in\mathcal N_{\Omega_i}(x^*), i=1,...,n$, where $a_{i,j}$ is the $(i,j)$th entry of the
adjacency matrix of $\mathcal G$. Hence, there exist $g(\mathbf{x}^*)\in\partial\mathbf{f}(\mathbf{x}^*)$ and {$\lambda^*\in\mathbb R^{nq}$} such that $-g(\mathbf{x}^*)-\alpha\mathbf{L\lambda}^*\in\mathcal N_{\Omega}(\mathbf{x}^*)$, equivalently,
$\mathbf{0}_{nq}=P_{\mathcal T_{\Omega} (\mathbf{x}^*)}(-g(\mathbf{x}^*)-\alpha\mathbf{L\lambda}^*)$.
(\ref{optimal_condition_1}) is proved.
\end{IEEEproof}

\subsection{Distributed Continuous-Time Projected Algorithm}

For the optimization problem (\ref{optimization_problem}), we propose a distributed optimization algorithm as follows:
\begin{subequations}\label{distributed feedback_v2_x}
\begin{align}
\dot x_i(t) &  =   P_{\mathcal T_{\Omega_i} (\mathbf{x}_i(t))}\bigg{[}- g_i(x_i(t))-\alpha\sum_{j=1}^{n}  a_{i,j}(x_i(t)-x_j(t))\nonumber\\
&\hspace{0 cm}-\alpha\sum_{j=1}^{n}a_{i,j}(\lambda_i(t)-\lambda_j(t))\bigg{]},
\, g_i(x_i(t))\in \partial{f^i}(x_i(t)),\\
\dot \lambda_i(t) &= \alpha \sum_{j=1}^{n}a_{i,j}(x_i(t)-x_j(t)), \label{distributed feedback_v2_v}
\end{align}
\end{subequations}
where $t\geq 0,\,i\in\{1,\ldots,n\},\,x_i(0)=x_{i0}\in\Omega_i\subset\mathbb R^q$, $\lambda_i(0)=\lambda_{i0}\in\mathbb R^q$, {$\alpha >0$}, and $a_{i,j}$ is the $(i,j)$th element of the 
adjacency matrix of graph $\mathcal G$, $\mathcal T_{\Omega_i}(x_i(t))$ is the tangent cone of $\Omega_i$ at an element $x_i(t)\in\Omega_i$ and $P_{\mathcal T_{\Omega_i}(x_i(t))}(\cdot)$ is the projection operator to $\mathcal T_{\Omega_i}(x_i(t))$.

\begin{remark}
Algorithm (\ref{distributed feedback_v2_x}) is motivated by the primal-dual type continuous-time algorithms, which was firstly proposed in \cite{WE:2011} and later on extended in \cite{GC:2014,YHL:SCL:2015,KCM:A:2015,LW:TAC:2015}.
If the state constraints are relaxed to $\Omega_i=\mathbb R^q,\,i\in\{1,\ldots,n\}$, then algorithm (\ref{distributed feedback_v2_x}) is consistent with the algorithm proposed in Section IV of \cite{GC:2014}. Algorithm (\ref{distributed feedback_v2_x}) also incorporates projection operation to handle constraints, which had also been adopted in \cite{QLX:A:2016} and \cite{LW:TAC:2015}.  However, \cite{QLX:A:2016} only handled homogeneous constraints, and \cite{LW:TAC:2015} may produce unbounded states, which may be hard to implement in practice. Here our proposed algorithm (\ref{distributed feedback_v2_x}) handles the problems with local constraints and can guarantee the boundedness of states.
\hfill$\blacklozenge$
\end{remark}

\section{Main results}\label{convergence_simulation}
In this section, we first introduce additional preliminaries for nonsmooth analysis, and then give the convergence analysis of the algorithm with an illustrative simulation.
\subsection{Nonsmooth Analysis}

To study our algorithm, we need concepts related to nonsmooth analysis.  Consider a differential inclusion \cite{AC:1984} in the form of
\begin{eqnarray}\label{DI}
  \dot x(t) \in  \mathcal H (x(t)),\quad x(0)=x_0,\quad t\geq 0,
\end{eqnarray}
where $\mathcal H:\mathbb R^q\rightarrow \mathfrak B(\mathbb R^q)$ is a set-valued map  with nonempty compact values.
Let $\tau>0$. A solution of \eqref{DI} defined on $[0, \tau] \subset [0,\infty)$ is an absolutely continuous function $x:[0, \tau]\rightarrow \mathbb R^q$ such that \eqref{DI} holds for almost all $t\in [0, \tau]$  ({in the sense of Lebesgue measure}).
Recall that the solution $t\mapsto x(t)$ to (\ref{DI}) is a \textit{right maximal solution} if it cannot be extended forward in time. We assume that all right maximal solutions to (\ref{DI}) exist on $[0,\infty)$. A set $\mathcal M$ is said to be {\em weakly invariant} \cite{R:1998} (resp., {\em strongly invariant}) with respect to (\ref{DI}) if  $\mathcal M$ contains a maximal solution \cite{R:1998} (resp., all maximal solutions) of (\ref{DI}) for every $x_0\in\mathcal M$.   A point $x_*$ is an {\em almost cluster point} \cite[p.~311]{AC:1984} of a measurable function $\phi(\cdot)$ when $t\rightarrow\infty$ if $\mu\{t\geq 0:\|\phi(t)-x_*\|\leq \varepsilon\}=\infty$ for all $\varepsilon>0$, where $\mu(\cdot)$ is the Lebesgue measure.

Let $\mathcal D$ be a compact, strongly positive invariant set with respect to (\ref{DI}). Let $W$ be a nonnegative {\em lower semicontinuous} (see \cite[p.~22]{AC:1984}) function defined on $\mathbb{R}^q\times\mathbb{R}^q$ 
and $V$ be a nonnegative lower semicontinuous and {\em inf-compact} (see \cite[p.~292]{AC:1984}) function defined on $\mathbb R^q$.
{
Assume there exists an upper semicontinuous ({
see \cite[p.~41]{AC:1984}}) map $ \tilde{\mathcal H}(x)$ with closed values such that $\mathcal H (x)\subset \tilde{\mathcal H}(x)$ for all $x\in\mathcal D$ and $\mathbf{0}_q\in\tilde{\mathcal H}(x)$ if and only if $\mathbf{0}_q\in{\mathcal H}(x)$,} we introduce a result  for the existence of an almost cluster point.
\begin{lemma}\label{almost_convergence}
If $ \phi(\cdot) \in\mathbb R^{q}$ is a solution of (\ref{DI}) with $ \phi(0) = x_0\in\mathcal D$ such that $$V(\phi(t))-V(\phi(s))+\int_{s}^{t}W(\phi(\tau),\dot \phi(\tau))\mathrm{d}\tau\leq 0,\quad t\geq s\geq 0,$$
then $\phi(\cdot)$ and $\dot \phi(\cdot)$ have almost cluster points $x_*$ and $v_*$, which satisfy  $W(x_*,v_*)=0$. If, in addition, $W(x,v)>0$ for all $x\in\mathbb R^q$ and all $v\not= \mathbf{0}_{q}$, then $x_*$ is an equilibrium of the differential inclusion (\ref{DI}).
\end{lemma}
\begin{IEEEproof}
{
By \cite[Proposition 5,~p.~311]{AC:1984}, $\phi(\cdot)$ and $\dot \phi(\cdot)$ have almost cluster points $x_*$ and $v_*$ which satisfy  $W(x_*,v_*)=0$.

If, in addition, $W(x,v)>0$ for all $x\in\mathbb R^q$ and all $v\not= \mathbf{0}_{q}$, then $v_*=\mathbf{0}_q$. Let $\{t_i\}_{i=1}^{\infty}$ be a increasing nonnegative sequence such that $t_i\rightarrow\infty$ and $\{\phi(t_i),\dot \phi(t_i)\}\rightarrow (x_*,\mathbf{0}_q)$. Clearly, $\dot \phi(t_i)\in\mathcal{H}(\phi(t_i))\subset\tilde{\mathcal H}(\phi(t_i))$ for all $i\in\{1,2,\ldots,\infty\}$. Because $\tilde{\mathcal H}(\cdot)$ is upper semicontinuous, $\mathbf{0}_{q}\in \tilde{\mathcal H}(x_*)$ by definition. Recall that $\mathbf{0}_{q}\in \tilde{\mathcal H}(x_*)$ is equivalent to $\mathbf{0}_{q}\in{\mathcal H}(x_*)$, $x_*$ is an equilibrium of the differential inclusion (\ref{DI}).}
\end{IEEEproof}

Furthermore, we introduce a lemma, which is inspired by \cite[Proposition 3.1]{HHB:TAC:2009} and is used in the convergence analysis.

\begin{lemma}\label{converge}
Let $\mathcal D$ be a compact, strongly positive invariant set with respect to (\ref{DI}), and $ \phi(\cdot) \in\mathbb R^{q}$ be a solution of (\ref{DI}) with $ \phi(0) = x_0\in\mathcal D$. If $z$ is an almost cluster point of $\phi(\cdot)$ and a Lyapunov stable equilibrium of (\ref{DI}), then $z= \lim_{t\rightarrow\infty} \phi(t)$.
\end{lemma}

\begin{IEEEproof}
Suppose $z$ is an almost cluster point of $\phi(\cdot)$ and $z$ is Lyapunov stable. Let $\varepsilon>0$.
Since $z$ is Lyapunov stable, there exists $\delta=\delta(\varepsilon, z)>0$ such that the solution $\tilde \phi(t)$ of system (\ref{DI}) with $\tilde \phi(0) = y\in\mathcal B_{\delta}(z)$ satisfies that $\tilde \phi(t)\in\mathcal B_{\varepsilon}(z)$ for all $t\geq 0$.
Since $z$ is an almost cluster point of $\phi(\cdot)$, there exists $h = h(\delta, x_0)>0$ such that $\phi(h)\in\mathcal B_{\delta}(z)$. It follows from our construction of $\delta $ that $\phi(t)\in\mathcal B_{\varepsilon}(z)$ for all $t\geq h$.
Because $\varepsilon>0$ is arbitrary, $z = \lim_{t\rightarrow\infty} \phi(t)$.
\end{IEEEproof}

\subsection{Convergence Analysis}
Let $\mathbf {x} \triangleq [x_1^{\rm T},\ldots,x_n^{\rm T}]^{\rm T}\in\Omega\subset\mathbb R^{nq}$ and $ \lambda \triangleq [\lambda_1^{\rm T},\ldots,\lambda_n^{\rm T}]^{\rm T}\in\mathbb R^{nq}$ with $\Omega\triangleq\prod_{i=1}^n \Omega_i$. Algorithm (\ref{distributed feedback_v2_x}) can be written in a compact form
\begin{eqnarray}
\begin{bmatrix}\dot {\mathbf x}(t) \\ \dot {\lambda}(t) \end{bmatrix}
 \in  \mathcal F(\mathbf{x}(t),\lambda (t)),\,\,\mathbf{x}(0)=\mathbf{x}_0\in\Omega,\,\, \lambda (0)=\lambda _0\in\mathbb R^{nq}, \label{feedback_comp}
\end{eqnarray}
where $ \mathcal F(\mathbf{x},\lambda )\triangleq \Bigg{\{}\begin{bmatrix} P_{\mathcal T_{\Omega} (\mathbf{x})}[- \alpha\mathbf{L}\mathbf {x} -  \alpha \mathbf{L}\lambda   -  g(\mathbf{x})]\\
\alpha \mathbf{L}\mathbf{x}\end{bmatrix}: g(\mathbf{ x})\in\partial \mathbf{f}(\mathbf {x} )\Bigg{\}}$ and $\mathbf{L}=L_n\otimes I_q\in\mathbb R^{nq\times nq}$.

\begin{remark}\label{Remark2}
The optimization algorithm (\ref{feedback_comp}) is of the form $\dot x(t)\in P_{\mathcal T_{K} ({x}(t))}[\mathcal H(x(t))]$, where $x(0)=x_0\in K$, $K$ is a closed convex subset of $\mathbb R^q$, and $\mathcal H$ is an {\em upper semicontinuous}  map with nonempty compact convex values. It follows from Proposition 2 of \cite[p.~266]{AC:1984} and Theorem 1 of \cite[p.~267]{AC:1984} that algorithm (\ref{feedback_comp}) has right maximal solutions on $[0, \infty)$.  Since $P_{\mathcal T_{K} ({x}(t))}[\mathcal H(x(t))]\subset \mathcal T_{K} ({x}(t))$, $K$ is a strongly invariant set to $\dot x(t)\in P_{\mathcal T_{K} ({x}(t))}[\mathcal H(x(t))]$.
{
In addition, $P_{\mathcal T_{K} ({x}(t))}[\mathcal H(x(t))]\subset \mathcal H(x(t)) -\mathcal N_{K}(x(t))$, $\mathbf{0}_q\in P_{\mathcal T_{K} ({x}(t))}[\mathcal H(x(t))]$ if and only if $\mathbf{0}_q\in \mathcal H(x(t)) -\mathcal N_{K}(x(t))$, and $\mathcal H(x(t)) -N_{K}(x(t))$ is upper semicontinuous because both $\mathcal H(x(t))$ and $\mathcal N_{K}(x(t))$ are upper semicontinuous. Hence, Lemma \ref{almost_convergence} can be applied to the convergence analysis of algorithm (\ref{feedback_comp}).
}
\hfill$\blacklozenge$
\end{remark}

Because $L_n$ is symmetric by Assumption \ref{Assumption1},
$L_n$ can be factored as $L_n=Q\Lambda Q^{\rm T}$ by the symmetric eigenvalue decomposition, where $Q$ is an orthogonal matrix and $\Lambda$ is a diagonal matrix whose diagonal entries are the eigenvalues of $L_n$.
Define a diagonal matrix $\overline \Lambda\in\mathbb{R}^{n\times n}$ such that $\overline \Lambda_{i,i}=1/\Lambda_{i,i}$ if $\Lambda_{i,i}>0$ and $\overline\Lambda_{i,i}=2k\alpha$ if $\Lambda_{i,i}=0$ for $i\in\{1,\ldots,n\}$.
The following lemma provides a result when {$\alpha>0$ and $0<k< \frac{1}{\alpha\lambda_{\max}(L_n)}$}.

\begin{lemma}\label{L_N_COM}
Consider algorithm \eqref{feedback_comp} under Assumption \ref{Assumption1} with $0<k< \frac{1}{\alpha\lambda_{\max}(L_n)}$. Then $Q_n=k\alpha^2 Q(\frac{1}{k\alpha}\overline \Lambda -I_n)Q^{\rm T}> 0$ and $\alpha L_n-k\alpha^2 L_n^2=L_n Q_n L_n$.
\end{lemma}

\begin{IEEEproof}
With $0<k< \frac{1}{\alpha\lambda_{\max}(L_n)}$, it is easy to prove $Q_n>0$.

Because $L_n=Q\Lambda Q^{\rm T}$ and $\Lambda\overline \Lambda \Lambda=\Lambda$ by the definition of $\overline \Lambda$,
\begin{eqnarray*}
L_nQ_nL_n &=& k\alpha^2 L_nQ(\frac{1}{k\alpha}\overline \Lambda -I_n)Q^{\rm T}L_n\\
&=& k\alpha^2 Q\Lambda Q^{\rm T} \big{[} Q(\frac{1}{k\alpha}\overline \Lambda -I_n)Q^{\rm T} \big{]} Q\Lambda Q^{\rm T}\\
&=& \alpha Q \Lambda\overline \Lambda \Lambda Q^{\rm T}-k\alpha^2( Q\Lambda Q^{\rm T})^2\\
&=&\alpha Q \Lambda Q^{\rm T}-k\alpha^2( Q\Lambda Q^{\rm T})^2\\
&=& \alpha L_n-k\alpha^2 L_n^2
\end{eqnarray*}
which implies the conclusion.
\end{IEEEproof}

If  3) of Assumption \ref{Assumption1} holds, there exists $(\mathbf{x}^*,\lambda^*)\in\Omega\times \mathbb{R}^{nq}$ satisfying (\ref{optimal_condition}) by Lemma \ref{optimal_eqv}. Let $\mathbf{x}^*\in\Omega$ and $\lambda^*\in\mathbb{R}^{nq}$ be the vectors such that (\ref{optimal_condition}) is satisfied.  Define
\begin{eqnarray}
V_1^*(\mathbf{x},\lambda) & \triangleq & \frac{1}{2}\|\mathbf{x}-\mathbf{x}^*\|^2+\frac{1}{2}\|\lambda-\lambda^*\|^2, \label{V}\\
V_2^*(\mathbf{x},\lambda)& \triangleq & \mathbf{f}(\mathbf{x})-\mathbf{f}(\mathbf{x}^{*})+ \alpha\frac{1}{2}\mathbf{x}^{\rm T}\mathbf{Lx}+\alpha\mathbf{x}^{\rm T}\mathbf{L}{\lambda}.\label{V_2}
\end{eqnarray}

\begin{remark}
{Functions $V_1^*(\mathbf{x},\lambda)$ and $V_2^*(\mathbf{x},\lambda)$ are constructed to form the candidates of Lyapunov functions in the theoretical analysis.
Function $V_1^*(\mathbf{x},\lambda)$ is also used as a Lyapunov function in \cite{GC:2014}  to prove algorithm convergence of unconstrained distributed optimization, which is a very good result.
In the analysis of \cite{GC:2014}, the   cost function was assumed to have a finite number of critical points and the quadratic Lyapunov functions were used.
However, in this note, the cost functions are assumed to be convex, which means that the cost function may have infinitely many solutions (or infinitely many critical points). Function $V_2^*(\mathbf{x},\lambda)$ uses the convexity  property to tackle convex cost functions (see part ($iii$) and ($iv$) of proof to Lemma \ref{differential_eqn}).
}\hfill$\blacklozenge$
\end{remark}

Recall that if $\phi(\cdot)$ is a solution of (\ref{DI}) and $V:\mathbb{R}^{q}\to\mathbb{R}$ is locally Lipschitz and {\em regular} (see \cite[p.~39]{Clarke:1983}), then $\dot \phi(t)$ and $\dot V(\phi(t))$ exist {\em almost everywhere}. Next, we give the following result, whose proof is given in Appendix.

\begin{lemma}\label{differential_eqn}
Suppose Assumption \ref{Assumption1} holds. Let $V_1^*(\mathbf{x},\lambda)$ and $V_2^*(\mathbf{x},\lambda)$ be as defined in \eqref{V} and \eqref{V_2}, and let $(\mathbf{x}(t),\lambda(t))$ be a  trajectory to  algorithm  (\ref{distributed feedback_v2_x}) or (\ref{feedback_comp}).
\begin{itemize}
\item[($i$)]
    $\dot V_1^*(\mathbf{x}(t),\lambda (t))\leq -\alpha \mathbf {x}^{\rm T} (t)\mathbf{L}\mathbf {x} (t)\leq 0$ for almost all $t\geq 0$.
\item[($ii$)]
    $\dot V_2^*(\mathbf{x}(t),\lambda (t))\leq -\|\dot{\mathbf {x}}(t)\|^{2} +\alpha^2 {\mathbf{x}}^{\rm T}(t) \mathbf{L}^2\mathbf{x}(t)$ for almost all $t\geq 0$.

\item[($iii$)] Let $0<k< \frac{1}{\alpha\lambda_{\max}(L_n)}$. The function $V^*(\mathbf{x},\lambda )=V_1^*(\mathbf{x},\lambda )+kV_2^*(\mathbf{x},\lambda )$ is nonnegative with all $(\mathbf{x},\lambda )\in \Omega\times \mathbb{R}^{nq}$.
\item[($iv$)]  With $V^*(\mathbf{x},\lambda )$ defined in part ($iii$) for $0<k< \frac{1}{\alpha\lambda_{\max}(L_n)}$, $\dot V^*(\mathbf{x}(t),\lambda (t))\leq -k\|\dot {\mathbf{x}}(t)\|^2-\dot {\lambda}^{\rm T}(t)\mathbf{Q}\dot {\lambda}(t)\leq 0$  for almost all $t\geq 0$, where $\mathbf{Q}\in\mathbb{R}^{nq\times nq}$ is positive definite.
\end{itemize}
\end{lemma}

Based on Lemmas \ref{converge} and \ref{differential_eqn}, we obtain our main result for state boundedness and convergence of the proposed algorithm.

\begin{theorem}\label{theorem_convergence}
Suppose Assumption \ref{Assumption1} holds and let $(\mathbf{x}(t),\lambda(t))$ be a trajectory to algorithm  (\ref{distributed feedback_v2_x}) or (\ref{feedback_comp}). Then
\begin{itemize}
\item[($i$)] $(\mathbf{x}(t),\lambda(t))$ is bounded;

\item[($ii$)] $(\mathbf{x}(t),\lambda(t))$ converges to a point $(\bar{\mathbf{x}},\bar \lambda)$ such that $\bar{\mathbf{x}}=\mathbf{1}_n \otimes\bar x$ and $\bar x$ is an optimal solution to problem \eqref{optimization_problem}.
\end{itemize}
\end{theorem}

\begin{IEEEproof}
In this theorem, part ($i$) claims that an equilibrium point of algorithm \eqref{feedback_comp} is Lyapunov stable and any trajectory of algorithm \eqref{feedback_comp} is bounded; part ($ii$) further claims that any  trajectory of algorithm \eqref{feedback_comp} converges to one of the equilibria of algorithm \eqref{feedback_comp}.

($i$) Let $V_1^*(\mathbf{x},\lambda)$ be as defined in \eqref{V}. It is clear that $V_1^*(\mathbf{x},\lambda)$ is positive definite, $V_1^*(\mathbf{x},\lambda)=0$ if and only if $(\mathbf{x},\lambda)=(\mathbf{x}^*,\lambda^*)$, and $V_1^*(\mathbf{x},\lambda)\rightarrow\infty $ as $(\mathbf{x},\lambda)\rightarrow\infty$.

By ($i$) of Lemma \ref{differential_eqn}, $\dot V_1^*(\mathbf{x}(t),\lambda (t)) \leq 0$  for almost all $t\geq 0$. Hence, $\mathcal D\triangleq\{(\mathbf{x},\lambda)\in\Omega\times\mathbb R^{nq}:V_1^*(\mathbf{x},\lambda)\leq M\}$, where $M>0$, is strongly positive invariant. Note that $V_1^*(\cdot,\cdot)$ is positive definite and {$V_1^*(\mathbf{x},\lambda)\rightarrow\infty $ as $(\mathbf{x},\lambda)\rightarrow\infty$}. Set $\mathcal D$ is bounded and the solution $(\mathbf{x}(t),\lambda(t))$ is also bounded. Part ($i$) is thus proved.

($ii$)  Let $V^*(\mathbf{x},\lambda)$ be as defined in ($iii$) of Lemma \ref{differential_eqn}. Due to ($iv$) of Lemma \ref{differential_eqn}, $\dot V^*(\mathbf{x}(t),\lambda (t))\leq -k\|\dot {\mathbf{x}}(t)\|^2-\dot {\lambda}^{\rm T}(t)\mathbf{Q}\dot {\lambda}(t)\leq 0$  for almost all $t\geq 0$, where $\mathbf{Q}\in\mathbb{R}^{nq\times nq}$ is positive definite.
Define $W(\dot {\mathbf{x}}, \dot \lambda)=k\|\dot {\mathbf{x}}\|^2+\dot {\lambda}^{\rm T}\mathbf{Q}\dot {\lambda}$. It is clear that $W(\dot {\mathbf{x}}, \dot \lambda)=0$ if and only if $\dot {\mathbf{x}}=\mathbf{0}_{nq}$ and $\dot \lambda=\mathbf{0}_{nq}$.

Recall that $(\mathbf{x}(t),\lambda(t))$ is bounded by ($i$) and $V^*(\mathbf{x},\lambda)$ is inf-compact and nonnegative with all $(\mathbf{x},\lambda )\in \Omega\times \mathbb{R}^{nq}$ by ($iii$) of Lemma \ref{differential_eqn}.
Note that
\begin{eqnarray*}
V^*(\mathbf{x}(t),\lambda(t))-V^*(\mathbf{x}(s),\lambda (s)) &=& \int_{s}^{t}\dot V^*(\mathbf{x}(\tau),\lambda(\tau))\mathrm{d}\tau\\
&\leq & -\int_{s}^{t} W(\dot {\mathbf{x}}(\tau), \dot \lambda(\tau))\mathrm{d}\tau.
\end{eqnarray*}
By Lemma \ref{almost_convergence}, $(\mathbf{x}(t),\lambda(t))$ has an almost cluster point $(\bar{\mathbf{x}},\bar \lambda)\in\Omega\times\mathbb R^{nq}$ and $(\bar{\mathbf{x}},\bar \lambda)$ is an equilibrium point of \eqref{feedback_comp}.

Define a function
$\bar V(\mathbf{x},\lambda) \triangleq \frac{1}{2}\|\mathbf{x}-\bar{\mathbf{x}}\|^2+\frac{1}{2}\|\lambda-\bar \lambda\|^2$.
It is clear that $\bar V(\mathbf{x},\lambda)$ is positive definite, $\bar V(\mathbf{x},\lambda)=0$ if and only if $(\mathbf{x},\lambda)=(\bar{\mathbf{x}},\bar \lambda)$, and $\bar V(\mathbf{x},\lambda)\rightarrow \infty$ if $(\mathbf{x},\lambda)\rightarrow\infty$. Because $(\bar{\mathbf{x}},\bar \lambda)$ is an equilibrium point of \eqref{feedback_comp}, $(\bar{\mathbf{x}},\bar \lambda)$ satisfies \eqref{optimal_condition}.
Moreover, it follows from ($i$) of Lemma \ref{differential_eqn} that $\bar V(\mathbf{x}(t),\lambda(t))$ along the trajectories of (\ref{distributed feedback_v2_x}) satisfies $\dot{\bar V}(\mathbf{x}(t),\lambda (t))\leq 0$  for almost all $t\geq 0$.
Hence, $(\bar{\mathbf{x}},\bar \lambda)$ is a Lyapunov stable equilibrium point to the system (\ref{distributed feedback_v2_x}).

Clearly, $(\bar{\mathbf{x}},\bar \lambda)$ is an almost cluster point of $(\mathbf{x}(t),\lambda(t))$ and $(\bar{\mathbf{x}},\bar \lambda)$ is a Lyapunov stable equilibrium. According to Lemma \ref{converge}, $(\mathbf{x}(t),\lambda(t))$ converges to $(\bar {\mathbf{x}},\bar \lambda)$ as $t\rightarrow\infty$.  Because $(\bar{\mathbf{x}},\bar \lambda)$ is  an equilibrium point of \eqref{feedback_comp}, there exists $\bar x\in\Omega_0\subset\mathbb R^q$ such that $\bar{\mathbf{x}}=\mathbf{1}_n \otimes\bar x$ and $\bar x$ is an optimal solution to problem \eqref{optimization_problem} by Lemma \ref{optimal_eqv}.
Part ($ii$) is thus proved.
\end{IEEEproof}

\begin{remark}
Theorem \ref{theorem_convergence} shows the convergence of the proposed algorithm.
The convergence analysis, in fact, can also be conducted following the method in \cite{SZH:CTT}.
\hfill$\blacklozenge$
\end{remark}

\begin{remark}
The convergence analysis in this note is based on {nonsmooth Lyapunov functions}, which can be regarded as an extension of the analysis on basis of smooth Lyapunov functions used in \cite{GC:2014,WE:2011,KCM:A:2015}.
Moreover, the novel technique proves that algorithm (\ref{distributed feedback_v2_x}) is able to solve optimization problems with a continuum of optimal solutions, and therefore, improves some previous ones in  \cite{WE:2011,KCM:A:2015}, which only handle problems with only one optimal point.
\hfill$\blacklozenge$
\end{remark}

\subsection{Numerical Simulation}
The following is a numerical example for illustration.

Consider the optimization problem (\ref{optimization_problem}) with $x\in\mathbb R$, where $\Omega_i=\{x\in\mathbb R:i-12\leq x\leq i-2\}$ and nonsmooth cost functions are
$$f^i(x) = \begin{cases}
             -x+i-5, & \mbox{if } x<i-5, \\
             0, & \mbox{if } i-5\leq x\leq i+5, \\
             x-i-5, & \mbox{if }x\geq i+5,
           \end{cases}\quad i=1,\ldots,5.$$
The adjacency matrix of the information sharing
graph $\mathcal G$ of algorithm (\ref{distributed feedback_v2_x}) is given by $$A = { \tiny\begin{bmatrix}
0 ~ 1 ~ 0 ~ 0 ~ 1\\
1 ~ 0 ~ 1 ~ 0 ~ 1\\
0 ~ 1 ~ 0 ~ 1 ~ 0\\
0 ~ 0 ~ 1 ~ 0 ~ 1\\
1 ~ 1 ~ 0 ~ 1 ~ 0
\end{bmatrix}}.$$ It can be easily verified that $\Omega_0=\cap_{i=1}^{5}\Omega_i=[-7,\,-1]$ and the optimal solution is $x=-1$, which is on the boundary of the constraint set $\Omega_0$. If there are no set constraints ($\Omega_i=\mathbb{R}$), every point in the set $[0,\,6]$ is an optimal solution.

\begin{figure}
  \centering
  \includegraphics[width=8 cm, height =6 cm]{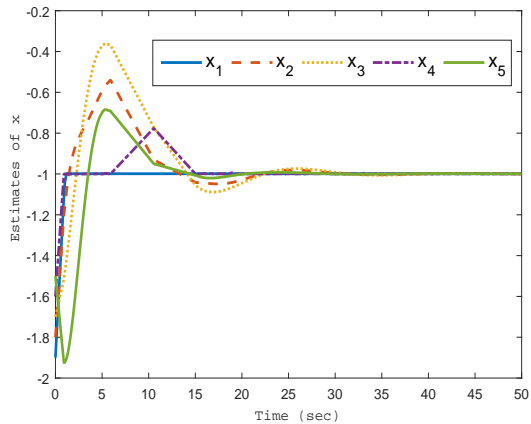}
  \caption{Trajectories of estimates for $x$ }
  \label{X}
\end{figure}

\begin{figure}
  \centering
  \includegraphics[width=8 cm, height =6 cm]{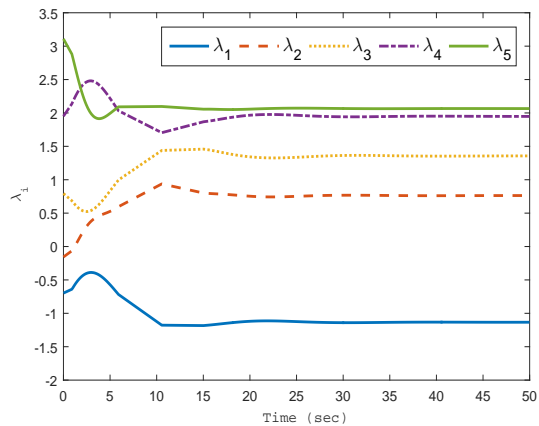}
  \caption{Trajectories of the auxiliary variable $\lambda$'s }
  \label{Lambda}
\end{figure}

\begin{figure}
  \centering
  \includegraphics[width=8 cm, height =6 cm]{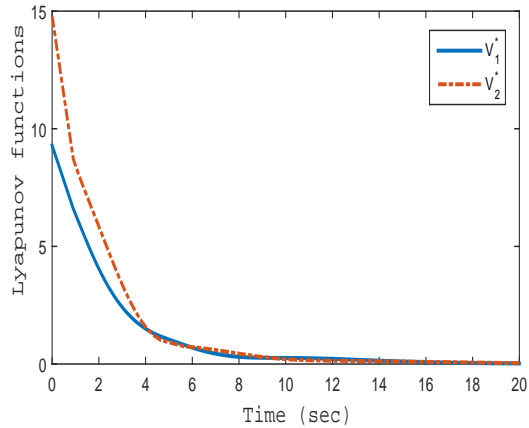}
  \caption{Trajectories of functions $V_1^*(\mathbf{x},\lambda)$ and $V_2^*(\mathbf{x},\lambda)$}
  \label{fig_funs}
\end{figure}

The trajectories of estimates for $x$ versus time are shown in Fig. \ref{X}.
It can be seen that all the agents converge to the same optimal solution
which satisfies all the local constraints and minimizes the sum of local cost functions,
without knowing other agents' constraints or feasible sets.
Fig. \ref{Lambda} shows the trajectories of the auxiliary variable $\lambda_i$'s and verifies the boundedness of the algorithm trajectories. {Fig. \ref{fig_funs} shows the trajectories of functions $V_1^*(\mathbf{x},\lambda)$ and $V_2^*(\mathbf{x},\lambda)$ versus time}.

\section{Conclusions}\label{conclusion}
In this note, a novel distributed projected continuous-time algorithm has been proposed for a distributed nonsmooth optimization under local set constraints.
By virtue of projected differential inclusions and nonsmooth analysis, the proposed algorithm has been proved to be convergent while keeping the states bounded.
Furthermore, based on the stability theory and convergence results for nonsmooth Lyapunov functions, the algorithm has been shown to solve the convex optimization problem with a continuum of optimal solutions.
Finally, the algorithm performance has also been illustrated via a numerical simulation.

\appendix[Proof of Lemma \ref{differential_eqn}]

($i$) Let $(\mathbf{x}(t),\lambda(t))$ be a trajectory to  algorithm  (\ref{distributed feedback_v2_x}) or (\ref{feedback_comp}). Recall that $\dot V_1^*(\mathbf{x}(t),\lambda (t))$ and $(\dot {\mathbf{x}}(t),\dot {\lambda} (t))$ exist for almost all $t\geq 0$. Suppose $\dot V_1^*(\mathbf{x}(t),\lambda (t))$ and $(\dot {\mathbf{x}}(t),\dot {\lambda} (t))$ exist at a positive time instant $t$. By \eqref{feedback_comp}, there exists $g(\mathbf{x}(t))\in\partial \mathbf{f}(\mathbf {x} (t))$ such that $\dot {\mathbf {x} }(t)=P_{\mathcal T_{\Omega} (\mathbf{x}(t))}[- \alpha\mathbf{L}\mathbf {x} (t)-  \alpha \mathbf{L}\lambda (t)  -  g(\mathbf{x}(t))]$ and $\dot{\lambda}(t)=\alpha \mathbf{L}\mathbf{x}(t)$.

Clearly, $\dot {\mathbf {x} }(t)=P_{\mathcal T_{\Omega} (\mathbf{x}(t))}[- \alpha\mathbf{L}\mathbf {x} (t)-  \alpha \mathbf{L}\lambda (t)  -  g(\mathbf{x}(t))]$ implies
$$- \alpha\mathbf{L}\mathbf {x} (t)-  \alpha \mathbf{L}\lambda (t)  -  g(\mathbf{x}(t))-\dot {\mathbf {x} }(t)\in \mathcal N_{\Omega}(\mathbf{x}(t)),$$
where $\mathcal N_{\Omega}(\mathbf{x}(t))\triangleq\{d\in\mathbb R^{nq}:d^{\rm T}(\tilde {\mathbf{x}}-\mathbf{x}(t))\leq 0,\,\forall \tilde {\mathbf{x}}\in\Omega \}$ is the normal cone of $\Omega$ at an element $\mathbf{x}(t)\in\Omega$. Hence,
\begin{eqnarray*}
\big{(} \alpha\mathbf{L}\mathbf {x} (t)+  \alpha \mathbf{L}\lambda (t)  +  g(\mathbf{x}(t))+\dot {\mathbf {x} }(t)\big{)}^{\rm T}\big{(}\mathbf{x}(t)-\tilde {\mathbf{x}}\big{)}\leq 0 ,
\end{eqnarray*}for all $\tilde {\mathbf{x}}\in\Omega$.

By choosing $\tilde {\mathbf{x}}= {\mathbf{x}^*}$,
\begin{eqnarray}\label{in_qn1}
\big{(} \alpha\mathbf{L}\mathbf {x} (t)+  \alpha \mathbf{L}\lambda (t)  +  g(\mathbf{x}(t))+\dot {\mathbf {x} }(t)\big{)}^{\rm T}\big{(}\mathbf{x}(t)- \mathbf{x}^*\big{)}\leq 0.
\end{eqnarray}

By Assumption \ref{Assumption1} and (\ref{optimal_condition_2}), we have $\mathbf{L}= \mathbf{L}^{\rm T}$ and $\mathbf{Lx}^{*}=\mathbf{0}_{nq}$, therefore,
\begin{eqnarray}\label{in_qn2}
\dot {\mathbf {x} }^{\rm T}(t)\big{(}\mathbf{x}(t)- \mathbf{x}^*\big{)} & \leq & -\alpha \mathbf {x}^{\rm T} (t)\mathbf{L}\mathbf {x} (t)-\alpha \mathbf {x}^{\rm T}(t)\mathbf{L}\lambda (t)\nonumber\\
&& -g(\mathbf{x}(t)) ^{\rm T}\big{(}\mathbf{x}(t)- \mathbf{x}^*\big{)}.
\end{eqnarray}

Furthermore, it follows from $\dot{\lambda}(t)=\alpha \mathbf{L}\mathbf{x}(t)$ that
\begin{eqnarray}\label{in_qn3}
\frac{1}{2}\frac{\mathrm d}{\mathrm{d}t}\|\lambda(t)-\lambda^*\|^2=\alpha(\lambda(t)-\lambda^*)^{\rm T} \mathbf{L}\mathbf{x}(t).
\end{eqnarray}

In view of \eqref{in_qn2} and \eqref{in_qn3},
\begin{eqnarray}\label{in_qn4}
\frac{\mathrm d}{\mathrm{d}t}V_1^*(\mathbf{x}(t),\lambda (t)) &\leq & -\alpha \mathbf {x}^{\rm T} (t)\mathbf{L}\mathbf {x} (t)-g(\mathbf{x}(t)) ^{\rm T}\big{(}\mathbf{x}(t)- \mathbf{x}^*\big{)}\nonumber\\
&&-\alpha {\lambda^*}^{\rm T} \mathbf{L}\mathbf{x}(t)\nonumber\\
&=& -\big{(}g(\mathbf{x}(t)) -g(\mathbf{x}^*)\big{)}^{\rm T}\big{(}\mathbf{x}(t)- \mathbf{x}^*\big{)}\nonumber\\
&&-(g(\mathbf{x}^*)+\alpha\mathbf{L\lambda}^*)^{\rm T}(\mathbf{x}(t)-\mathbf{x}^*)\nonumber\\
&&-\alpha \mathbf {x}^{\rm T} (t)\mathbf{L}\mathbf {x} (t),
\end{eqnarray}
where $g(\mathbf{x}^*)\in\partial\mathbf{f}(\mathbf{x}^*)$ is chosen such that $P_{\mathcal T_{\Omega} (\mathbf{x}^*)}(-g(\mathbf{x}^*)-\alpha\mathbf{L\lambda}^*)=\mathbf{0}_{nq}$.

Note that $P_{\mathcal T_{\Omega} (\mathbf{x}^*)}(-g(\mathbf{x}^*)-\alpha\mathbf{L\lambda}^*)=\mathbf{0}_{nq}$ implies  $-g(\mathbf{x}^*)-\alpha\mathbf{L\lambda}^*\in \mathcal N_{\Omega}(\mathbf{x}^*)$, where $\mathcal N_{\Omega}(\mathbf{x}^*)$ is the normal cone of $\Omega$ at an element $\mathbf{x}^*\in\Omega$. Hence,
$$(-g(\mathbf{x}^*)-\alpha\mathbf{L\lambda}^*)^{\rm T}({p}-\mathbf{x}^*)\leq 0$$
for all ${p}\in \Omega$.
Since $\mathbf{x}(t)\in\Omega$, we have
\begin{eqnarray}\label{in_qn5}
(-g(\mathbf{x}^*)-\alpha\mathbf{L\lambda}^*)^{\rm T}(\mathbf{x}(t)-\mathbf{x}^*)\leq 0.
\end{eqnarray}
Because $\mathbf{f}(\mathbf{x})$ is convex, $\big{(}g(\mathbf{x}(t)) -g(\mathbf{x}^*)\big{)}^{\rm T}\big{(}\mathbf{x}(t)- \mathbf{x}^*\big{)}\geq 0$ with $g(\mathbf{x}(t))\in\partial\mathbf{f}(\mathbf{x}(t))$ and $g(\mathbf{x}^*)\in\partial\mathbf{f}(\mathbf{x}^*)$. It follows from \eqref{in_qn4} that
\begin{eqnarray}
\frac{\mathrm d}{\mathrm{d}t}V_1^*(\mathbf{x}(t),\lambda (t))\leq -\alpha \mathbf {x}^{\rm T} (t)\mathbf{L}\mathbf {x} (t)\leq 0.
\end{eqnarray}

($ii$) Let $(\mathbf{x}(t),\lambda(t))$ be a trajectory  to  algorithm  (\ref{distributed feedback_v2_x}) or (\ref{feedback_comp}). Recall that $\dot V_2^*(\mathbf{x}(t),\lambda (t))$ and $(\dot {\mathbf{x}}(t),\dot {\lambda} (t))$ exist for almost all $t\geq 0$.  Suppose $\dot V_2^*(\mathbf{x}(t),\lambda (t))$ and $(\dot {\mathbf{x}}(t),\dot {\lambda} (t))$ exist at a positive time instant $t$. Since $\mathbf{f}(\mathbf{x})$ is convex in $\mathbf{x}$,
\begin{eqnarray*}
\mathbf{f}(\mathbf{x}(t))-\mathbf{f}(\mathbf{x}(t-h)) &\leq & \langle p, \mathbf{x}(t)- \mathbf{x}(t-h)\rangle,\\
\mathbf{f}(\mathbf{x}(t+h))-\mathbf{f}(\mathbf{x}(t)) &\geq & \langle p, \mathbf{x}(t+h)- \mathbf{x}(t)\rangle.
\end{eqnarray*}
for all $p\in\partial \mathbf{f}(\mathbf{x}(t))$ and $h\in(0,t]$.

Dividing both sides of the inequalities by $h\in(0,t]$ and letting  $h\rightarrow 0$, we obtain
\begin{eqnarray}\label{dot_phi_1}
 \frac{\mathrm{d}}{\mathrm{d}t} \mathbf{f}(\mathbf{x}(t)) = \langle p, \dot {\mathbf{x}}(t)\rangle,\quad \forall p\in\partial \mathbf{f}(\mathbf{x}(t)).
\end{eqnarray}

By \eqref{feedback_comp}, there exists $g(\mathbf{x}(t))\in\partial  \mathbf{f}(\mathbf{x}(t))$ such that $\dot {\mathbf {x} }(t)=P_{\mathcal T_{\Omega} (\mathbf{x}(t))}[- \alpha\mathbf{L}\mathbf {x} (t)-  \alpha \mathbf{L}\lambda (t)  -  g(\mathbf{x}(t))]$ and $\dot{\lambda}(t)=\alpha\mathbf{Lx}(t)$. Choose $p=g(\mathbf{x}(t))$. Then $\frac{\mathrm{d}}{\mathrm{d}t} \mathbf{f}(\mathbf{x}(t))=g(\mathbf{x}(t))^{\rm T}\dot {\mathbf{x}}(t).$
Hence,
\begin{eqnarray} \label{Dif_V_2}
\frac{\mathrm d}{\mathrm{d}t}V_2^*(\mathbf{x}(t),\lambda (t)) &=& [ \alpha\mathbf{L}\mathbf {x} (t)+  \alpha \mathbf{L}\lambda (t)  +  g(\mathbf{x}(t))]^{\rm T}\dot {\mathbf{x}}(t)\nonumber\\
&&+\alpha^2\mathbf{x}^{\rm T}(t)\mathbf{L}^2 \mathbf{x}(t).
\end{eqnarray}

Set $K = \mathcal T_{\Omega} (\mathbf{x}(t))$, $v=\mathbf{0}_{nq}\in K$, $ u= -[ \alpha\mathbf{L}\mathbf {x} (t)+  \alpha \mathbf{L}\lambda (t)  +  g(\mathbf{x}(t))]\in\mathbb R^{nq}$, and $P_{K}(u)=\dot {\mathbf{x}}(t)$ in (\ref{projection_inequality}).
It follows from (\ref{projection_inequality}) that
$$[ \alpha\mathbf{L}\mathbf {x} (t)+  \alpha \mathbf{L}\lambda (t)  +  g(\mathbf{x}(t))]^{\rm T}\dot {\mathbf{x}}(t)\leq -\|\dot {\mathbf{x}}(t)\|^2.$$
Hence, $\frac{\mathrm d}{\mathrm{d}t}V_2^*(\mathbf{x}(t),\lambda (t)) \leq -\|\dot {\mathbf{x}}(t)\|^2+\alpha^2\mathbf{x}^{\rm T}(t)\mathbf{L}^2 \mathbf{x}(t)$, which follows from \eqref{Dif_V_2}.

($iii$) Let $0<k< \frac{1}{\alpha\lambda_{\max}(L_n)}$ and note that $\mathbf{Lx}^*=\mathbf{L}^{\rm T}\mathbf{x}^* =\mathbf{0}_{nq}$. It can be easily verified that
\begin{eqnarray*}
V^*(\mathbf{x},\lambda ) &=& V_1^*(\mathbf{x},\lambda )+k V_2^*(\mathbf{x},\lambda )\\
&=& J_1(\mathbf{x},\lambda)+J_2(\mathbf{x})+J_3(\mathbf{x}),
\end{eqnarray*}
where $J_1(\mathbf{x},\lambda)=\frac{1}{2}\|\mathbf{x}-\mathbf{x}^*\|^2+\frac{1}{2}\|\lambda-\lambda^*\|^2 + k\alpha(\mathbf{x}-\mathbf{x}^*)^{\rm T}\mathbf{L}(\lambda-\lambda^*),$ $J_2(\mathbf{x})= k\alpha\frac{1}{2}\mathbf{x}^{\rm T}\mathbf{Lx},$ and $J_3(\mathbf{x})=k[\mathbf{f}(\mathbf{x})-\mathbf{f}(\mathbf{x}^{*}) +\alpha(\mathbf{x}-\mathbf{x}^*)^{\rm T}\mathbf{L}{\lambda^*}].$
To prove $V^*(\mathbf{x},\lambda )$ is nonnegative for all $(\mathbf{x},\lambda )\in \Omega\times \mathbb{R}^{nq}$, we show $J_1(\mathbf{x},\lambda)\geq 0$, $J_2(\mathbf{x})\geq 0$, and $J_3(\mathbf{x})\geq 0$ for all $(\mathbf{x},\lambda )\in \Omega\times \mathbb{R}^{nq}$.

Since $\mathbf{L}$ is positive semi-definite, 
\begin{eqnarray}\label{J_2}
J_2(\mathbf{x})= k\alpha\frac{1}{2}\mathbf{x}^{\rm T}\mathbf{Lx}\geq 0,
\end{eqnarray}
and $((\mathbf{x}-\mathbf{x}^*)+(\lambda-\lambda^*))^{\rm T}\mathbf{L}((\mathbf{x}-\mathbf{x}^*)+(\lambda-\lambda^*))\geq 0$ for all $(\mathbf{x},\lambda )\in \Omega\times \mathbb{R}^{nq}$.
Hence,
\begin{eqnarray}\label{semi_posi_ineq}
&&(\mathbf{x}-\mathbf{x}^*)^{\rm T}\mathbf{L}(\mathbf{x}-\mathbf{x}^*)+(\lambda-\lambda^*)^{\rm T}\mathbf{L}(\lambda-\lambda^*)\geq \nonumber  \\
&&\hspace{2 cm}-(\mathbf{x}-\mathbf{x}^*)^{\rm T}(\mathbf{L}+\mathbf{L}^{\rm T})(\lambda-\lambda^*).
\end{eqnarray}
Let $\mu_i,\,i=1,\ldots,n,$ be the eigenvalues of $L_n\in\mathbb R^{n\times n}$. Since the eigenvalues of $I_q$ are $1$, it follows from the properties of Kronecker product that the eigenvalues of $\mathbf{L}=L_n\otimes I_q$ are $\mu_i\times 1,\,i=1,\ldots,n$. Thus, $\lambda_{\max}(L_n)=\lambda_{\max}(\mathbf{L})$.

Because of Assumption \ref{Assumption1}, $\mathbf{L}=\mathbf{L}^{\rm T}$. By \eqref{semi_posi_ineq},
\begin{eqnarray*}
k\alpha(\mathbf{x}-\mathbf{x}^*)^{\rm T}\mathbf{L}(\lambda-\lambda^*) &\geq & -\frac{k\alpha}{2} (\mathbf{x}-\mathbf{x}^*)^{\rm T} \mathbf{L}(\mathbf{x}-\mathbf{x}^*)\\
&&-\frac{k\alpha}{2} (\lambda-\lambda^*)^{\rm T}\mathbf{L}(\lambda-\lambda^*)\\
&\geq &  -\frac{k\alpha\lambda_{\max}(L_n)}{2}  \|\mathbf{x}-\mathbf{x}^*\|^2  \\ &&-\frac{k\alpha\lambda_{\max}(L_n)}{2}\|\lambda-\lambda^*\|^2.
\end{eqnarray*}
Due to  $1-k\alpha\lambda_{\max}(L_n)>0$,
\begin{eqnarray}\label{J_1}
J_1(\mathbf{x},\lambda) & \geq & \frac{1}{2}(1-k\alpha\lambda_{\max}(L_n))\|\mathbf{x}-\mathbf{x}^*\|^2\nonumber\\
&&+\frac{1}{2}(1-k\alpha\lambda_{\max}(L_n))\|\lambda-\lambda^*\|^2\geq 0.
\end{eqnarray}

Since $\mathbf{f}(\mathbf{x})$ is convex in $\mathbf{x}\in\Omega$,
\begin{eqnarray*}
J_3(\mathbf{x}) &=& k[\mathbf{f}(\mathbf{x})-\mathbf{f}(\mathbf{x}^{*}) +\alpha(\mathbf{x}-\mathbf{x}^*)^{\rm T}\mathbf{L}{\lambda^*}]\\
  & \geq & k[ (p+\alpha\mathbf{L}{\lambda^*})^{\rm T}   (\mathbf{x}-\mathbf{x}^*) ],\quad \forall p\in\partial \mathbf{f}(\mathbf{x}^*).
\end{eqnarray*}

Note that there exists $g(\mathbf{x}^*)\in\partial\mathbf{f}(\mathbf{x}^*)$ such that $P_{\mathcal T_{\Omega} (\mathbf{x}^*)}(-g(\mathbf{x}^*)-\alpha\mathbf{L\lambda}^*)=\mathbf{0}_{nq}$, which follows from \eqref{optimal_condition_1}. Choose $p\triangleq g(\mathbf{x}^*)$. In light of \eqref{in_qn5} and similar arguments above \eqref{in_qn5},  $$(p+\alpha\mathbf{L}{\lambda^*})^{\rm T}   (\mathbf{x}-\mathbf{x}^*)\geq 0$$
for all $\mathbf{x}\in\Omega$ with $p\triangleq g(\mathbf{x}^*)$. Hence,
\begin{eqnarray}\label{J_3}
J_3(\mathbf{x}) \geq 0,\quad \forall \mathbf{x}\in\Omega.
\end{eqnarray}

In view of \eqref{J_2}, \eqref{J_1}, and \eqref{J_3},   $V^*(\mathbf{x},\lambda )=V_1^*(\mathbf{x},\lambda )+kV_2^*(\mathbf{x},\lambda )$ is nonnegative with all $(\mathbf{x},\lambda )\in \Omega\times \mathbb{R}^{nq}$.

($iv$) It follows from part ($i$) and ($ii$) that $\dot V^*(\mathbf{x}(t),\lambda (t))\leq - \mathbf {x}^{\rm T} (t)[\alpha\mathbf{L}-k\alpha^2\mathbf{L}^2]\mathbf {x} (t) -k\|\dot{\mathbf {x}}(t)\|^{2} $ for almost all $t\geq 0$.

{With $Q_n>0$  as defined in  Lemma \ref{L_N_COM}, we have $$L_nQ_nL_n=\alpha L_n-k\alpha^2 L_n^2$$
by  Lemma \ref{L_N_COM}. Define $\mathbf{Q} = Q_n\otimes I_q>0$.
Recalling $\dot \lambda (t)=\alpha \mathbf{L}\mathbf{x}(t)$, {it can be easily proved  that  $$\mathbf {x}^{\rm T}(t)(\alpha\mathbf{L}-k\alpha^2\mathbf{L}^2)\mathbf {x}(t)=\mathbf {x}^{\rm T}(t)\mathbf{L}\mathbf{Q}\mathbf{L}\mathbf {x}(t)=\dot {\lambda}^{\rm T}(t)\mathbf{Q}\dot {\lambda}(t).$$}Hence, $\dot V^*(\mathbf{x}(t),\lambda (t))\leq -k\|\dot {\mathbf{x}}(t)\|^2-\dot {\lambda}^{\rm T}(t)\mathbf{Q}\dot {\lambda}(t)\leq 0$  for almost all $t\geq 0$.}
\bibliographystyle{IEEEtran}

%




\end{document}